\documentclass[11pt]{amsart}
\usepackage{amsfonts,amssymb}
\usepackage{epsfig}
\usepackage{latexsym}
\usepackage[all]{xy}
\usepackage{amsmath}
\usepackage{amsfonts}
\usepackage{amssymb}

\def\R{\mathbb R}

\def\Z{\mathbb Z}

\def\S{\mathbb S}

\def\Q{\mathbb Q}

\def\V{\mathcal V}
\def\W{\mathcal W}
\def\A{\mathcal A}
\newtheorem{Def}{Definition}[section]
\newtheorem{Thm}{Theorem}[section]
\newtheorem{Cor}{Corollary}[section]

\newtheorem{Prop}{Proposition}[section]

\newtheorem{Remark}{Remark}[section]
\begin{document}
\title[Index of a finitistic space]{Index of a finitistic space and a generalization of the topological central point theorem}
\author{Satya Deo*}
\thanks{* While this work was done, the author was supported by the DST (Govt of India) Research Grant sanction letter number SR/S4/MS:567/09 dated 18.02.2010}
\date{}
\address{Harish-Chandra Research Institute, \\
Chhatnag Road, Jhusi,\\
Allahabad 211 019, India. }
\email{sdeo@mri.ernet.in,  vcsdeo@yahoo.com }
\subjclass[2000]{}
\begin{abstract}
 In this paper we prove that if $G$ is a $p$-torus (resp. torus) group acting without fixed points on a finitistic space $X$ (resp. with finitely many orbit types), then the $G$-index $i_G(X) < \infty$. Using this G-index we obtain a generalization of the Central Point Theorem and also of the Tverberg Theorem for any d-dimensional Hausdorff space.

\end{abstract}
\vskip 1in
\maketitle

\section{Introduction}
Let $X$ be a topological space with an involution T and let $\S^n$ denote the $n$-sphere with antipodal involution. Then the $\Z_2-$index of the $\Z_2-$space $X$ is defined as the smallest integer n such that there is an equivariant map $X\to \S^n$ (see \cite{Mat}, p.95). Then the Borsuk-Ulam theorem asserts that the $\Z_2-$index of $\S^n$ is $n$. The idea of this $\Z_2-$index has been generalized in many ways, first for any finite group (see \cite{con}) and then for arbitrary $G-$spaces where $G$ is a compact Lie group acting on a space $X$ (see \cite{Mat}, \cite{vol}). Its multiple applications in proving several nonembedding results like the van Kampen-Flores theorem, theorems like Central Point theorem and the Tverberg theorem etc. make it a powerful tool in topological combinatorics. The idea of an "ideal-index" of a G-space $X$, using the spectral sequence arising from the principal G-bundle map of Borel construction, has proved to be quite effective (see \cite {vol}, \cite{kar}). This ideal-index defines very naturally several numerical indices, in particular the numerical G-index $i_G(X)$ of a G-space $X$. If $f:X\to Y $ is a G-map then it follows easily that $i_G(X)\leq i_G(Y).$ This result, combined with the G-index of the deleted join with natural group actions, is utilized in showing the nonexistence of a G-map from a G-space $X$ to a G-space $Y$. However, in doing so it is necessary to know that the G-index of $Y$ is finite. The first aim of this paper is to show that if $G$ is a $p$-torus (resp. torus) group acting without fixed points on a finitistic space $X$ (resp. with finitely many orbit types), then the $G$-index $i_G(X) < \infty$. This is done in Section 2.  Then, as an application, we obtain in Section 5, a generalization of the Central Point Theorem and also of the Tverberg Theorem for any d-dimensional Hausdorff space. The cohomology groups of the classifying spaces of these groups have a neat presentation with a polynomial part in several variables which can be used as a multiplicative set in understanding the scalar product.

\section{Preliminaries}
Let $G$ be a compact Lie group acting continuously on a topological space $X$. Suppose $B_G$ denotes the classifying space of $G$ and $E_G\to B_G$ is the principal $G$-bundle of $G$. We have two continuous maps $\pi_1,\;\pi_2$ defined by the Borel construction
$$\begin{array}{ccccc}
 X& \leftarrow & X\times E_G & \rightarrow & E_G\\
\downarrow & & \downarrow & & \downarrow\\
X\backslash G & \stackrel{\pi_2}\leftarrow & X\times_G E_G & \stackrel{\pi_1}\rightarrow & B_G
\end{array}$$
Let $H^*$ denote the Alexander-Spanier cohomology and $H^{*}_{G}(X)=H^*(X\times_GE_G),\;H^{*}_{G}(pt)=H^{*}(B_G)$ denote the corresponding equivariant cohomology of the $G$-space $X$ and the trivial $G$-space $\{pt\}$ consisting of a point. Then $\Lambda^*=H^{*}_G(pt)$ is a commutative ring with identity and $H^{*}_{G}(X)$ becomes a $\Lambda^*$-module via the ring homeomorphism induced by the equivariant map $X\to\{pt\}$. The Serre spectral sequence of the bundle map $\pi_1:X_G\to B_G$ with fibre $X$ defines a sequence of homeomorphism
$$\Lambda^*\to E^{*0}_{2}\to\cdots\to E^{*0}_{r}\to\cdots\to E^{*0}_{\infty}\subset H^{*}_{G}(X)$$
The kernel of the composite map $\Lambda^*\to E^{*0}_{\infty}$ or $\Lambda^*\to H^{*}_{G}(X)$, denoted by $\hbox{Ind}^{G}(X)$, is called the ``ideal index'' of the $G$-space $X$. The kernel of the homeomorphism $\Lambda^*\to E_{r+1}$, denoted by $_r\hbox{Ind}^{G}(X)$, is called the ${r}^{th}$ filtration of the ideal index $\hbox{Ind}^G(X)= _{\infty}\hbox{Ind}^G(X)$. Now we have (see \cite{vol})

{\Def  Let $\alpha (\neq 0)\in\hbox{Ind}^G(X)$ and $i_X(\alpha)$ denote the smallest $r$ such that $\alpha\in\; _r\hbox{Ind}^G(X)$. Then $\min \{r\;|\; i_X(\alpha)=r$ and $\alpha (\neq 0)\in\hbox{Ind}^G(X)\}$ is called the $G$-index of the $G$ space $X$ and is denoted by $i_G(X)$. Similarly, $\min \{r\;|\; i_X(\alpha)=r$ and $\alpha (\neq 0)\in\hbox{Ind}^G(X)$ and is not a zero-divisor in $\Lambda^*\}$ is yet another $G$ index of $X$ and is denoted by $i'_G(X)$. }

\bigskip  Recall that the Yang's homology index, $\hbox{ind}_G(X)$, defined by using the equivariant homology theory with $\Z_2-$ coefficients is a well-known numerical $\Z_2-$index \cite{yan}. Clearly $i_G(X)\leq i'_G(X)$. If Ind$^G(X)=0$, we put $i_G(X)=\infty$. Note that if $X$ has a fixed point, then there is a monomorphism from $\Lambda^*\to H^{*}_{G}(X)$ and hence $i_G(X)=\infty$. If $f:X\to Y$ is a $G$-map, then $i_G(X)\leq i_G(Y)$. Also, if $G$ acts on the sphere $\S^k$ without fixed point, then $i_G(\S^k)=k+1$. It is obviously interesting to ask as to what are those $G$-spaces $X$ on which $G$ acts without fixed points and  for which $i_G(X)$ is finite. It must be pointed out that the index $i_G(X)$ of a $G$-space $X$ is essentially dependent on the action of $G$ on $X$, e.g., when the circle group $G=\S^1$ acts on an $n$-sphere $\S^n$ with fixed points, then $i_G(X)=\infty$, but when it acts without fixed points then $i_G(S^n)=n+1.$

 Here we show that if the group $G$  is a $p$-torus  $\Z^{r}_{p}$ or a torus $T^r$ (acting with finitely many orbit types) and $X$ is a finitistic $G$-space, then with suitable coefficient groups for the cohomology, $i_G(X)$ and $i'_G(X)$ both are finite. The result is already known when $X$  is a compact $G$-space or $X$ is a paracompact $G$-space of finite cohomological  dimension and $G$ is acting with finitely many orbit types (see \cite{vol}, Section 3).

\medskip Let us recall (see ~\cite{Br}, \cite{deo-tri}) that a paracompact Hausdorff space $X$ is said to be {\bf finitistic} if each open cover of $X$ has a finite dimensional open refinement. Every compact space is trivially finitistic and every paracompact space of finite covering dimension is also clearly finitistic. On the other hand, there are numerous finitistic spaces which are neither compact nor of finite cohomological dimension, e.g., the product space $X=(\prod^{\infty}_{n=1}\S^n)\times\R^k$. Finitistic spaces are the most general class of spaces suitable for cohomology theory of topological transformation groups (\cite{All},\cite{Br}). Hence it is natural to ask whether or not the index $i_G(X)$ is finite when $G$ is acting on a finitistic space $X$ without fixed points.

\section{Index of a finitistic $G$-space is finite}
Suppose the group $G=\Z^{r}_{p}$ or $T^r$ and the coefficients are in the fields $k=\Z_p$ or $\Q$ respectively. Then we know that (see \cite{Hsi}, p.45)

\begin{eqnarray*}
 H^*(B_G,k) & = & k[t_1,t_2,\cdots ,t_r],\; \deg t_i=2 (\hbox{resp.}\;1),\\
& & \;\hbox{when}\; G=T^r (\hbox{resp.}\; \Z^{r}_{2}),\; k=\Q (\hbox{resp.}\; \Z_2).\\
H^*(B_G,k) & = & k[t_1,t_2,\cdots ,t_r]\otimes \Lambda[v_1,\cdots, v_r],\; \deg v_i=1,\;t_i=\beta (v_i)\\
& & \hbox{when}\;G=\Z^{r}_{p}, \Q=\Z_p\;\hbox{and $p$ is odd}.
\end{eqnarray*}

Suppose $G$ acts on a finitistic space $X$ without fixed points. For each $x\in X$, the inclusion map $ G(x)\to X$ induces a continous map  $j_x: G(x)_G\to X_G$ in the Borel constructions. Let $j^{*}_{x}: H^{*}_{G}(X)\to H^{*}_{G}(G(x))$ denote the induced homomorphism in equivariant cohomology. Also, the constant $G$-map $G(x)\to \{pt\}$ induces a homomorphism $H^*(B_G)\to H^*(G(x))=H^*(B_G)$ for each $x\in X$. Let $S=R-\{0\}$ be the multiplicative subset where $R$ is the polynomial part of the ring $H^*(B_G, k)$. If $G$ acts on $X$ without fixed points, then $G_x$ is a proper subtorus of $G$ for each $x\in X$, and so from the cohomology algebra mentioned above it follows that for each $x$, there is an $s\in S$ which goes to zero under the map $H^*(B_G)\to H^*(B_{G_x})$. Note that $H^*(B_{G_x})$ depends only on the conjugacy class of $G_x$ in $G$.

\bigskip If $G=\Z^{r}_{p}$, then trivially there are only finitely many conjugacy classes of $G_x$, but if $G=T^r$, we {\bf assume} that there are only finitely many orbit types in $X$, and hence in this case also there are only finitely many conjugacy classes of $G_x$ in $G$. Choosing one $s_i\in S$ for each conjugacy class $G_{x_i},\;i=1,2,\cdots ,n$, we find that the product $s=s_1.s_2 \cdots s_n\in S$ has the property that $s$ goes to zero under the map $H^*(B_G)\to H^*(B_{G_x}),\;\;\forall\; x\in X$. If $\pi_{1}^{*}(s)\in H^{*}_{G}(X)$, we will show that $(\pi^{*}_{1}(s))^N=0$ for some integer $N$, which will prove that $i_G(X)<\infty$.

\bigskip Let us call the element $j^{*}_{x}(\pi^{*}_{1}(s))$ as the {\bf restriction} of the element $s$ on the orbit $G(x)$ and denote it as $s|G(x)$. Since $s|G(x)=0$ we find by the tautness property of the Alexander-Spanier cohomology that there is an invariant open neighborhood $U(x)$ of $x$ in $X$ such that $s|U(x)=0$. Thus we have now an open covering $\{U(x)\}$ of the space $X$. Since $X$ is finitistic, we find from \cite{deo-tri} that the orbit space $X/G$ is also finitistic. Hence the open covering $\{q(U(x))\}$ has a finite dimensional open refinement, say $\V=\{V_{\alpha}\;|\;\alpha\in\mathcal{A}\}$. Here $q:X\to X/G$ is the quotient map. Let $\{f_{\alpha}\;|\;\alpha\in\mathcal{A}\}$ be a partition of unity subordinate to $\V$ and $\W=\{W_{\alpha}\;|\;W_{\alpha}=f^{-1}_{\alpha}(0,\; 1]\subset V_{\alpha}\}$

\bigskip Then $\W$ is a finite-dimensional locally finite refinement of $\V$. Let $\{X_{\alpha}\;|\;\alpha\in\A\}$ be a shrinking of $\{W_{\alpha}\}$ which means $\{X_{\alpha}\;|\;\alpha\in\A\}$ is a locally finite, finite-dimensional covering of $X/G$ such that $\overline{X}_{\alpha}\subset W_{\alpha}$ for each $\alpha\in\A$. Put
$$Z_{\alpha}=q^{-1}(X_{\alpha}) \;\hbox{and}\;Y_{\alpha}=q^{-1}(W_{\alpha}).$$
Then $\{Z_{\alpha}\}$ and $\{Y_{\alpha}\}$ both are finite-dimensional locally finite coverings of the space $X$ consisting of invariant open sets so that $\{Z_{\alpha}\}$ is a shrinking of $\{Y_{\alpha}\}$

\bigskip Now for each $\alpha\in\A$, let $\nu_{\alpha}:X / G\to [0,1]$ be a continuous function ($X/G$ is normal) such that $\nu_{\alpha}(\overline{Z}_{\alpha})=1$ and $\nu_{\alpha}(X/G-Y_{\alpha})=0$. Let $\mu_{\alpha}$ be the composite map $\nu_{\alpha}\circ q: X\to [0,1]$. Then $\mu_{\alpha}(\overline{X}_{\alpha})=1$ and  $\mu_{\alpha}(X-W_{\alpha})=0.$
 Now for any finite nonempty subset $F$ of $\A_{\alpha}$, let
$$U(F)=\{x\in X\;|\;\min_{\alpha\in F}\mu_{\alpha}(x)>\max_{\alpha\notin F}\mu_{\alpha}(x)\}$$
Then $U(F)\subseteq\bigcap_{\alpha\in F}W_{\alpha}$,  and because of  finite dimensionality of $\{V_{\alpha}\}$, there is an integer $n$ such that $U(F)=\phi$ whenever $|F|\geq n$. Since $\{V_{\alpha}\}$ is locally finite,  and $\{W_{\alpha}\}$ is a covering,  $\{V(F)\;|\;F\subset\A\;\hbox{and}\;|F|<\infty\}$ is an invariant open cover of $X$. Moreover $|F_1|=|F_2|$ and $F_1\neq F_2$ implies $U(F_1)\bigcap U(F_2)=\phi$ because otherwise for some $x\in U(F_1)\bigcap U(F_2)$, we will have
\begin{eqnarray*}
 \min_{\alpha\in F_1}\mu_{\alpha}(x) & > & \max_{x\notin F_1}\mu_{\alpha}(x)\\
 &\geq & \min_{\alpha\in F_2}\mu_{\alpha}(x)\\
&\geq & \max_{\alpha\notin F_2}\mu_{\alpha}(x)\\
&\geq & \min_{x\in F_1}\mu_{\alpha}(x),
\end{eqnarray*}
a contradiction.

\bigskip  Let $U_i=\bigcup\{U(F)\;|\;|F|=i\}$. Then $\{U_1, U_2,\cdots , U_n\}$ is a covering of $X$ by invariant open sets. Since, for each $i$, $U_i$ is a disjoint union of $U(F_i)$, we have
$$H^{*}_{G}(U_i)=\Pi_{|F|=i}H^{*}_{G}(U(F))$$
Also, since $U(F)\subset\bigcap_{\alpha\in F}V_{\alpha}\subset V_{\alpha}$ for some $\alpha$, which is contained in $U(x)$, we find that $s|F=0$ for each $F$,  and therefore $s|U_i=0\;\;\forall\; i=1,2,\cdots ,n$. Considering the pair $(X,U_i)$ we observe that $\pi^{*}_{1}(s)\in H^{*}_{G}(X, U_i)$ and hence the product $\pi^{*}_{1}(s^n)= (\pi^{*}_{1}(s))^n\in H^{*}_{G}(X, X)=0$. Therefore $(\pi^{*}_{1}(s))^n$ is zero in $H^{*}_{G}(X)$. We have thus proved the following.

{\Thm Suppose $G=\Z^{r}_{p}$ (resp $T^r$) is a $p$-torus (resp. torus) acts on a finitistic space $X$ without fixed points. If $G=T^r$, we assume further that $G$ acts on $X$ with finitely many orbit types. Then with coefficients in $\Z_p$(resp. $\Q$), the index $i_G(X)$ of $X$ is finite.}

\bigskip The following result follows easily from definitions (see ~\cite{vol}).

{\Prop Let $X$ and $Y$ be two $G$-spaces.
\begin{enumerate}
 \item[(i)] If $\tilde{H}_i(X)=0$ for $i< N$, then $i_G(X)\geq N+1$.
\item[(ii)] If $H^i(Y)=0$ for $i>N-1$ and $i_G(Y)<\infty$, then $i_G(Y)\leq N$.
\end{enumerate}}

\bigskip We already know that if there is a $G$-map $f:X\to Y$, the $i_G(X)\leq i_G(Y)$. Hence we have the following

{\Cor Let $G$ be a $p$-torus acting on spaces $X,\;Y$ without fixed points. Suppose $\tilde{H}^i(X,\Z_p)=0$ for $i<N$ and $Y$ is a finitistic space such that $H^i(Y,\Z_p)=0$ for $i>N-1$. Then there is no $G$ map from $X$ to $Y$.}

{\Cor Let $G$ be a torus acting on spaces $X$ and $Y$ without fixed points. Assume that $H^i(X,\Q)=0$ for $i< N$ and $Y$ is finitistic having only finitely many orbit types such that $H^i(Y,\Q)=0$ for $i>N-1$. Then there is no $G$-map from $X$ to $Y$.}

{\Remark The condition that the space $X$ is finitistic in the above theorem is needed. For example take any group G and let $E_G\to B_G$ be the principal G-bundle with $E_G$ as the bundle space. Then $G$ acts on $E_G$ without fixed points (freely) and so the number of orbit types is finite. Since the space $(E_G)_G$ is paracompact and the fibres of $p:(E_G)_G\to B_G$ are acyclic, the Vietoris-Begle theorem implies that $p^*:H^*(B_G)\to H^*(E_G)$ is an isomorphism, i.e., $i_G(E_G)= \infty$.} The space $E_G$ is clearly nonfinitistic.

{\Remark The above theorem is true not only for the groups $\Z_p^r$ and $T^r$, but also for any compact connected Lie group $G$ provided the coefficients are taken in a field $k$ of characteristic zero, e.g., $k=\Q.$ In this case
$$ {\rm ker}  p^* = \bigcap_{x\in X} {\rm ker}  j_x^*p^*=\bigcap_{i=1}^{m}{\rm ker}  j_{x_i}^*p^*,    $$
where the number of orbit types is finite because $\alpha \in {\rm ker} j_x^*p^*$ for each $x\in X$ means $p^*(\alpha)\in {\rm ker} j_x^*$ for each $x\in X,$ and hence there is an integer n such that $(p^*(\alpha))^n =0,$ i.e., $\alpha^n\in {\rm ker}  p^*.$ But ${\rm ker}  p^*$ is a prime ideal implies $\alpha \in {\rm ker} p^*.$

Now if ${\rm ker}  p^*=0,$ then $\cap_{i=1}^m {\rm ker} j_{x_i}^*=0$  Hence there is an index k such that ${\rm ker} j_{x_k}p^* =0.$ This means the map $H^*(B_G)\to H^*(B_{G_{x_k}})$ is injective map, i.e., $G_{x_k}=G,$ which means $x_k$ is a fixed point, a contradiction. Hence $i_G(X)<\infty$.         }

\section{A generalization of the topological central point theorem}

The classical central point theorem says that ``{\it given any point set $X$ having $m=(d+1)(r-1)+1$ points in $\R^d$, there exists a point $x$ in $\R^d$ which is contained in the convex hull of every subset $F$ of $X$ having at least $d(r-1)+1$ points.}'' The point $x$ is called a central point of $X$. This theorem has a nontrivial  generalization known as the Tverberg Theorem which says that ``{\it given any finite set $X$ in $\R^d$  having $m=(d+1)(r-1)+1$ points, we can partition the set $X$ into $r$ subsets $X_1, X_2,\cdots , X_r$ such that} $$\bigcap^{r}_{i=1}\hbox{conv}(X_i)\neq \phi.\hbox{''}$$
Both of these theorems have been generalized in various directions (~\cite{Br}~\cite{Sar}, etc). The following is a recent generalization of the central point theorem~\cite{kar}.
{\Thm {\bf (Topological Central Point Theorem)} \\Let $m=(d+1)(r-1)$ and $\Delta^m$ be the $m$-simplex. Suppose $W$ is any $d$-dimensional metric space and $f:\Delta^m\to W$ is any continuous map then
$$\bigcap_{\dim F=d(r-1)}f(F)\neq \phi$$
where $F$ runs over all faces of $\Delta^m$ of dimension $d(r-1)$.}

\bigskip The following topological Tverberg theorem is also proved in the same paper ~\cite{kar} for the case when $r$ is a prime power. For general $r$ the theorem is still open. In contrast to this, the topological central point theorem is true for any $r$.

{\Thm {\bf (Topological Tverberg Theorem)}\\ Let $m=(d+1)r-1$ where $r$ is a prime power and let $\Delta^m$ be a $m$-dimensional simplex. Suppose $f:\Delta^m\to W$ is any continuous map to a $d$-dimensional metric space $W$. Then there exists $r$ disjoint faces $F_1,\cdots , F_r$ of $\Delta^m$ such that $$\bigcap^{r}_{i=1}f(F_i)\neq \phi .$$}

In what follows we will generalize both of the above theorems for continuous maps $f:\Delta^m\to W$ where $W$ is any Hausdorff space of dimension $d$, not necessarily a metric space.

\section{Generalization to $d$-dimensional Hausdorff spaces}
Suppose $\dim X$ denotes the covering dimension of a space $X$~\cite{Pea}. We know that if $A$ is a closed subspace of $X$, then $\dim A\leq \dim X$. If $P$ is any $m$-dimensional polyhedron then $\dim P=m$. We also point out that both of the above theorems are proved in \cite{kar} with full detail when the space $W$ is a polyhedron of dimension $d$, and the proof uses the $G$-index of a $G$-space $X$ as discussed earlier in this paper for $G=\Z_2$ or $G=(\Z_p)^{r}$.

\bigskip Let us first mention that the topological central point theorem  is not valid when the space $W$ is an arbitrary $d$-dimensional space. As an easy example, let us take wedge (one point union with $p$ as a common point) $W=\Delta^d\vee\Delta^m$ where $\Delta^d$ is the usual $d$-simplex and $\Delta^m$ is the usual $m$-simplex with indiscrete topology on $\Delta^m-\{p\}$. Note that $\Delta^m$ with indiscrete topology is  $0-$dimensional but not Hausdorff. Then the central point theorem cannot be true for the continuous map $f:\Delta^m\to W$ where $f:\Delta^m\to\Delta^d\vee\Delta^m$ is the obvious inclusion map into the the second part, because the intersection of all $d(r-1)$ dimensional faces of $\Delta^m$ is clearly empty.

\bigskip {We have now}

{\Thm Let $f:\Delta^m\to W$ be a continuous map where $m=(d+1)(r-1)$ and $W$ be a $d$-dimensional Hausdorff space. If the topological central point theorem is true for a $d$-dimensional polyhedron, then the theorem is true for any Hausdorff space $W$ of dimension $d$.}

\begin{proof}
 Since the map $f:\Delta^m\to W$ is continuous and $W$ is Hausdorff, $f(\Delta^m)$ is a compact closed subspace of $W$. Hence $\dim f(\Delta^m)\leq d$. Moreover, since $f:\Delta^m\to f(\Delta^m)$ is a perfect map, viz., the inverse image $f^{-1}(p)$ is compact for each $p$ in $p(W)$, the space $f(\Delta^m)$ is also metrizable (see~\cite{Pea} p. 96, Cor 5.8). Thus the space $M=f(\Delta^m)$ is a compact metric space of dimension at most $d$. Let us denote the map $f$ as the composite
$$\Delta^m\stackrel{f'}\to M\stackrel{i}\to W$$
 Since $M$ is a compact metric space of dimension at most $d$, given any $\epsilon > 0$, we can find a polyhedron $P_{\epsilon}$ of dimension at most $d$ and a map $g_{\epsilon}: M\to P_{\epsilon}$ (see~\cite{Eng}, p 322, 6.P.51) such that the diam($ g^{-1}_{\epsilon}(p))$ $\leq \epsilon$ for each $p\in P_{\epsilon}$.

\begin{picture}(300,135)(20,-40)
\put(70,30){\makebox(0,0)[t]{$\Delta^m$}}
\put(105,40){\makebox(0,0)[t]{$f'$}}
\put(145,30){\makebox(0,0)[t]{$M$}}
\put(85,25){\vector(1,0){50}}
\put(160,30){\vector(3,1){50}} \put(160,20){\vector(3,-1){50}}
\put(185,40){\makebox(0,0)[b]{$i$}}
\put(185,8){\makebox(0,0)[t]{$g_{\epsilon}$}}
\put(220,50){\makebox(0,0)[t]{$W$}}
\put(220,4){\makebox(0,0)[t]{$P_{\epsilon}$}}
\end{picture}

Now we know that the theorem is true for the map $g_{\epsilon}\circ f':\Delta^m\to P_{\epsilon}$ for any given $\epsilon$. This means $\exists\; p\in P_{\epsilon}$ such that
$$ p\in\bigcap g_{\epsilon}\circ f'(F),$$
  where $F$ runs over all faces of $\Delta^m$ of dimension $(d+1)r$. This implies $\cap f'(F)\neq\phi$. To see this suppose $F_1, F_2,\cdots ,F_k$ are all the $(d+1)r$ dimensional faces of $\Delta^m$. Suppose $f'(F_1),\cdots ,f'(F_n),\;n<k$ intersect and there exists an image set, $f'(F_{n+1})$ which is disjoint from $f'(F_1)\cap\cdots\cap f'(F_n)$. Let $$\delta=d(f'(F_{n+1}),\; \bigcap^{n}_{i=1}f'(F_i))>0.$$
Choose $\epsilon<\frac{\delta}{2}$, a polyhedron $P_{\epsilon}$ and a continuous map $g_{\epsilon}: M\to P_{\epsilon}$ such that diam ($ g^{-1}_{\epsilon}(p)$) $<\bar{\epsilon}, \; \forall\;p\in P_{\bar{\epsilon}}$. Since $\cap_{i\leq k}g_{\epsilon}\circ f'(F_i)\neq \phi$ we find that $\cap_{i\leq n+1}g_{\epsilon}\circ f(F_i)\neq \phi$, and so there is a point $x_1\in \cap^{n}_{i=1}f'(F_i)$ and a point $x_2\in f'(F_{n+1})$ such that $g(x_1)=g(x_2)=p\in P_{\epsilon}$. But this is a contradiction to the fact that diam($ g^{-1}_{\epsilon}(p)$) $<\epsilon$. Hence $\cap g_{\epsilon}\circ f'(F)\neq \phi$ and this means $\cap f(F)\neq\phi$.
\end{proof}

By a similar argument one can prove the Topological Tverberg Theorem also. Since the topological central point theorem and the topological Tverberg theorem both are already known for polyhedra ~\cite{kar}, we have the following:

{\Cor {\bf (The Topological Central Point Theorem for Hausdorff spaces)}:
Let $m=(d+1)(r-1)$ and $\Delta^m$ be a $m$-simplex. Suppose $W$ is a $d$-dimensional Hausdorff space and $f:\Delta^m\to W$ is a continuous map. Then if $F$ runs over all the faces of $\Delta^m$ of dimension $d(r-1)$, then
$$\bigcap_{\dim F=d(r-1)} f(F)\neq\phi.$$}

{\Cor {\bf (The Topological Tverberg Theorem for Hausdorff Spaces)}:
Let $m=(d+1)r-1$ where $r$ is a prime power and $\Delta^m$ be a $m$-simplex. Let $f:\Delta^m\to W$ be a continuous map into a $d$-dimensional Hausdorff space $W$. Then there exist $r$ disjoint faces $F_1,\cdots , F_r$ of $\Delta^m$ such that
$$\bigcap^{r}_{i=1}f(F_i)\neq\phi.$$}

\bigskip\noindent {\bf Acknowledgement : } The author is thankful to R. N. Karasev for a few e-mail discussions on his paper~\cite{kar} and also to the referee for some useful remarks.

\end{document}